\documentclass[11pt]{article}
\usepackage[utf8]{inputenc}
\usepackage[T2A]{fontenc}
\usepackage{amsmath,amsfonts,multicol,hyperref}
\usepackage{graphicx}
\usepackage{setspace} %
\usepackage{amssymb,amsthm,euscript}

\newcommand{\rr}{\mathbb{R}}
\newcommand {\nn} {\mathbb{N}}
\newcommand {\zz} {\mathbb{Z}}

 \newcommand {\al} {\alpha}
\newcommand {\be} {\beta}
\newcommand {\da} {\delta}

\newcommand {\la} {\lambda}

\newcommand {\sa} {\sigma}

\newcommand {\fy} {\varphi}

\newcommand{\IN}{{\subset}}

\newcommand {\mmm}{{\setminus}}
\newcommand{\8}{{\infty}}

\newcommand{\0}{{\varnothing}}
\newcommand{\vse}{$\blacksquare$}

\newcommand{\eK}{{\EuScript K}}
\newcommand{\mD}{{\mathcal D}}

\newtheorem{thm}{\bf Theorem}
 \newtheorem{cor}[thm]{\bf Corollary}
 \newtheorem{lem}[thm]{\bf Lemma}

\newcommand{\dok}{{\bf{Proof:  }}}

\title{On the connected components of  fractal cubes.}

\author{Dmitry Drozdov, Andrei Tetenov\\Gorno-Altaisk State University}

\begin{document}

\maketitle

\begin{abstract}
We show that  a fractal cube $F$ in $\mathbb R^3$ may have an uncountable set $Q$ of connected components $K_\alpha$ neither of which is contained in any plane, whereas the set $Q$ is a totally disconnected self-similar subset of the hyperspace $C(\mathbb R^3)$, isomorphic to a Cantor set. 
\end{abstract}

\smallskip
{\it 2010 Mathematics Subject Classification}. Primary: 28A80. \\
{\it Keywords and phrases.} fractal square, fractal cube, self-similar set, hyperspace

\bigskip

\subsection{Introduction.}

Take $n\geq 2$ and some subset  $\mD \IN \{0,1, \ldots , n-1\}^k$, $2\le\#\mD<n^k$ which we call  a digit set. Then there is unique non-empty compact set $F\IN \rr^k$ satisfying $$F=\frac{F+\mD}{n},$$ which we call fractal $k$-cube ( or fractal square if $k=2$).\\

Let $H=F+\zz^k$ and $H^c=\rr^k\mmm H$. We also denote $I=[0,1]^n$, $T(A):=\dfrac{A+\mD}{n}$ and $F_1=T(I)$.

For the case $k=2$ when $F$ is a fractal square, K.S.Lau, J.J.Luo and R.Hui \cite{LLR} proved  the following
  
\begin{thm}\label{LLR}
Let $F$ be a fractal square as in (1.1). Then $F$ satisfies either\\
(i) $H^c$ has a bounded component, which is also equivalent to: $F$ contains a non-trivial component that is not a line segment; or\\
(ii) $H^c$ has an unbounded component, then $F$ is either totally disconnected or all non-trivial components of $F$ are parallel line segments.
\end{thm}

We show that in  the case of fractal cubes in $\rr^3$, the situation is completely different. In our short note, we prove the following

\begin{thm}
There is a fractal cube  $F\IN \rr^3$ such  that $H^c$ is connected and $H$ is an uncountable union of unbounded components, each being invariant with respect to $\zz^3$-translations.
\end{thm}

We construct $F$ as a fractal cube with $n=5$ and a digit set $\mD\IN\{0,1,...,4\}^3$ which is a disjoint union of its subsets $\mD_0$ and $\mD_1$. These digit sets define disjoint connected fractal cubes $K_0$ and $K_1$.

The Hutchinson operators $T_i(A):=\dfrac{\mD_i+A}{5}$ for $\mD_i$  can be considered as  contraction maps $\widetilde{T_i}:C(\rr^3)\to C(\rr^3)$ whose Lipschitz constant is $1/5$. 
For any finite $\bar{\al}=\al_1\ldots\al_k\in\{0,1\}^k, k\in\nn$, we define $T_{\bar{\al}}=T_{\al_1}\circ\ldots\circ T_{\al_k}$ and $F_{\bar{\al}}=T_{\bar{\al}}\circ I.$ Let $\eK_{\bar{\al}}$ be the non-empty compact set satisfying $\eK_{\bar{\al}}=T_{\bar{\al}}(\eK_{\bar{\al}}).$ For any infinite string $\al=\al_1\al_2....\in \{0,1\}^\8$ we write $|\bar\al|=k$ and define \[K_\al=\bigcap\limits_{k=1}^\8T_{\al_1...\al_k}(I)\] In our case, each $K_\al$ is a  connected component of $F$ and 
\[F=\bigcup\limits_{\al\in\{0,1\}^\8}K_\al\].

Since all components $K_\al$ of $F$ are not line segments, the equivalence (i) of Theorem\ref{LLR}  does not hold. From the other side, the set $H^c$ is connected and unbounded, but all components of $F$ are not line segments, so (ii) does not work too.\\

Thus, the set $Q$ of connected components $\eK_\al$ of $F$ is uncountable and it is totally disconnected  if we consider it as a subset of the hyperspace $C(\rr^3)$.  

\begin{thm}
The set $Q\IN C(\rr^3)$ of connected components $\eK_\al$ of $F$ is a self-similar set generated by two contractions $\tilde{T_0}$ and $\tilde{T_1}$ of the hyperspace $C(\rr^3)$. There is a H\"older homeomorphism $\fy:Q\to C_{1/3}$ of the set $Q$ to
 the middle-third Cantor set $C_{1/3}$ which induces the isomorphism of self-similar structures on these sets.
\end{thm}

Let $\al=\al_1 \al_2 \ldots\in\{0,1\}^\8.$ For any $m \in\nn$ let $\la_{m}$ be the density of zeros in the word $\al_{1}\ldots\al_{m}.$ Put $\bar{\la}_\al=\limsup\limits_{m\to\8}{\la_{m}},\ \underline{\la}_\al=\liminf\limits_{m\to\8}{\la_{m}}$ and if $\underline{\la}_\al=\bar\la_\al$, we write $\la_\al=\lim\limits_{m\to\8}{\la_{m}}.$ If $\al$ is preperiodic with period length $p$, the $\la_\al$ is equal to the density of zeros in any period of $\al$.

We prove the following estimates for the dimension of the components $K_\al$:

\begin{thm}
For any $\al\in\{0,1\}^\8$ 
$$\underline{\dim}_B(\eK_\al)=\bar\la_\al\log_513+(1-\bar\la_\al)\log_544,$$
$$\bar{\dim}_B(\eK_\al)=\underline{\la}_\al\log_513+(1-\underline{\la}_\al)\log_544.$$
If $\al$ is preperiodic, then $$\dim_H(\eK_\al)=\dim_B(\eK_\al)=\la_\al \log_513+(1-\la_\al)\log_544.$$
\end{thm}

\subsection{Construction of the set $F$}

\includegraphics[width=.45\textwidth]{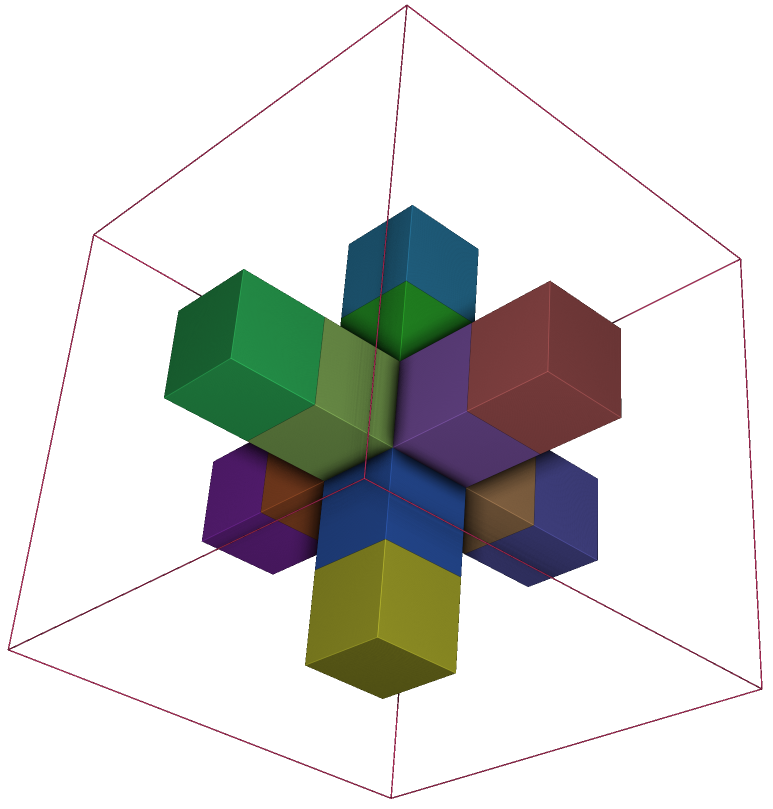}
\includegraphics[width=.45\textwidth]{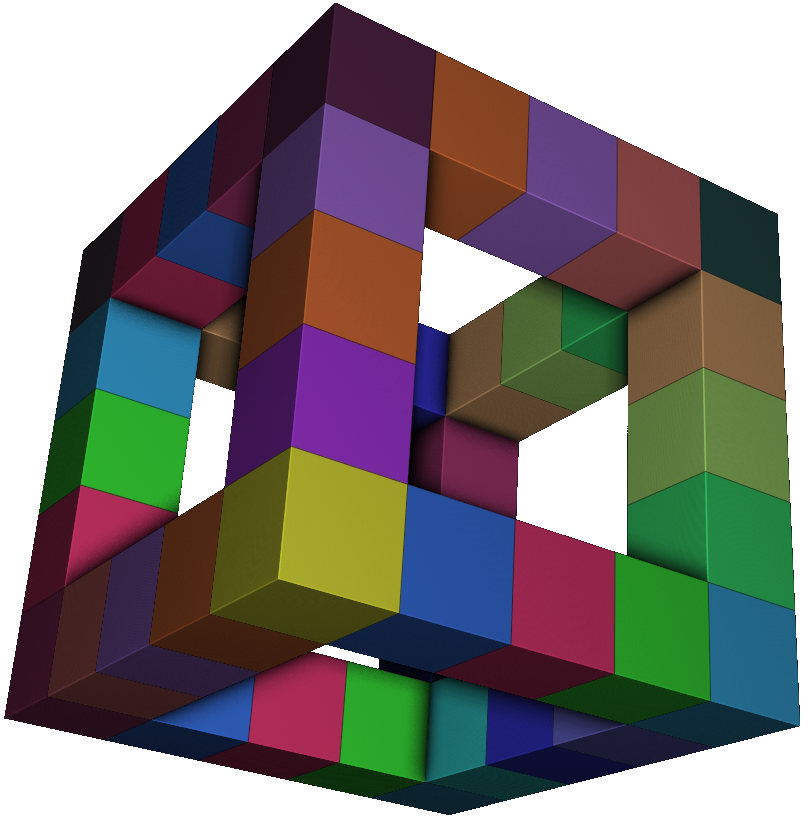}
\begin{center}
{\scriptsize Cross(left)  and  Frame (right).}    
\end{center}

Let $\mD_0$ and $\mD_1$ be the sets of the coordinates defining the cubes $\bar x+I$ forming subsets "Cross" and "Frame" of a cube $5\cdot I$, so that $\mD_0+I$ is a "Cross" and $\mD_1+I$ is a "Frame".

Let $T_0(A):=\dfrac{\mD_0+A}{5}$ and $T_1(A):=\dfrac{\mD_1+A}{5}$ be Hutchinson operators defined by $\mD_0$ and $\mD_1$. Let $\mD=\mD_0\bigcup\mD_1$ and $T(A):=\dfrac{\mD+A}{5}$. Let $\eK_0$ and $\eK_1$ be the fractal cubes corresponding to $\mD_0$ and $\mD_1$.\\

\includegraphics[width=.45\textwidth]{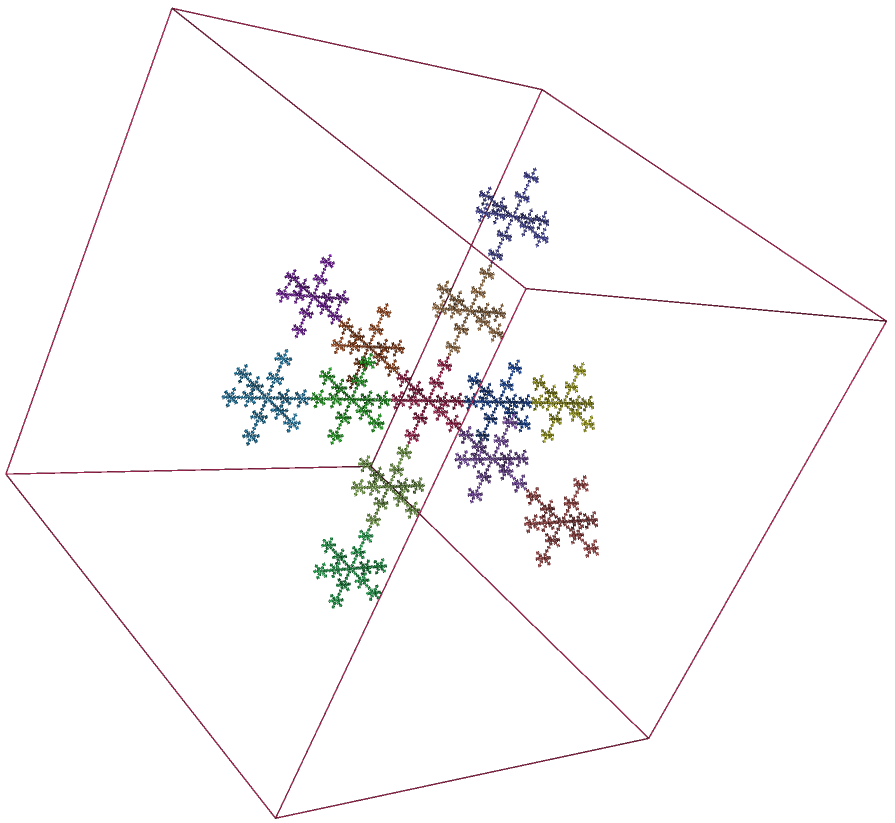}
\includegraphics[width=.45\textwidth]{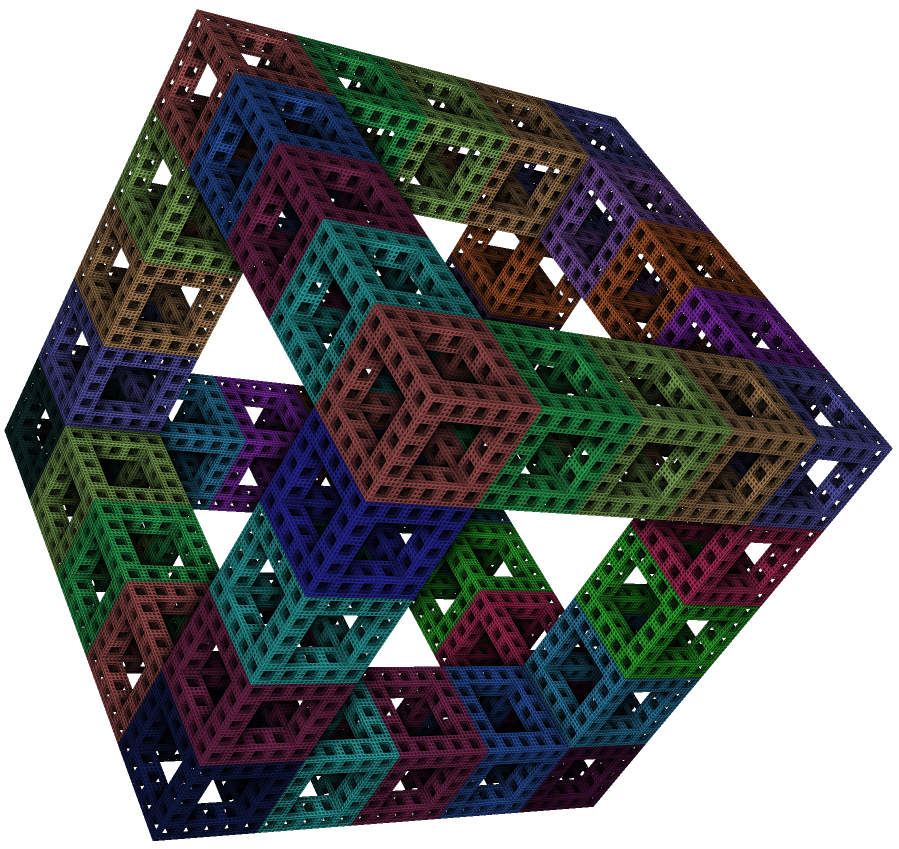}
\begin{center}
{\scriptsize Fractal cubes $\eK_0$(left)  and  $\eK_1$ (right).}    
\end{center}

Applying the approach, developed by L.Cristea and B.Steinsky in \cite{CS1}, we see  that since the set $\mD_0+I$ has exactly one pair of entrance and exit points on each pair of opposite faces of the cube $5I$ and has empty intersection with its edges and the intersection graph of  $\mD_0+I$ is a tree,  the fractal cube $\eK_0$ is a dendrite, all of whose ramification points have order 6.

The fractal cube $\eK_1$ is a $5\times 5\times 5$ version of  Menger sponge.\\

\includegraphics[width=.45\textwidth]{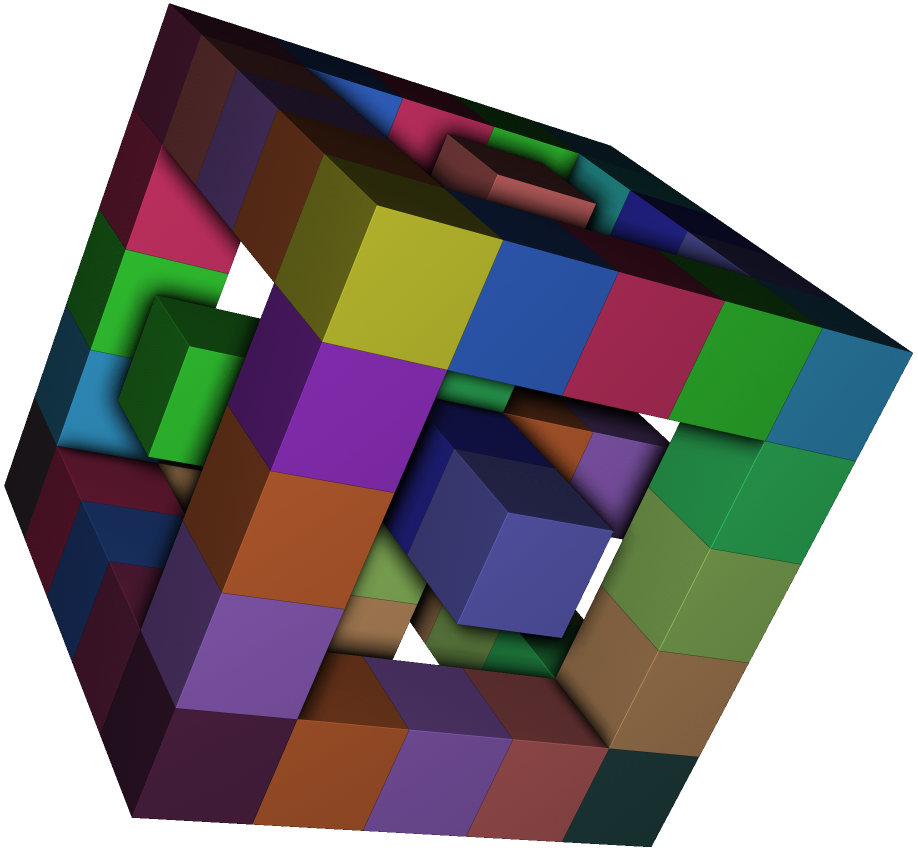}
\includegraphics[width=.45\textwidth]{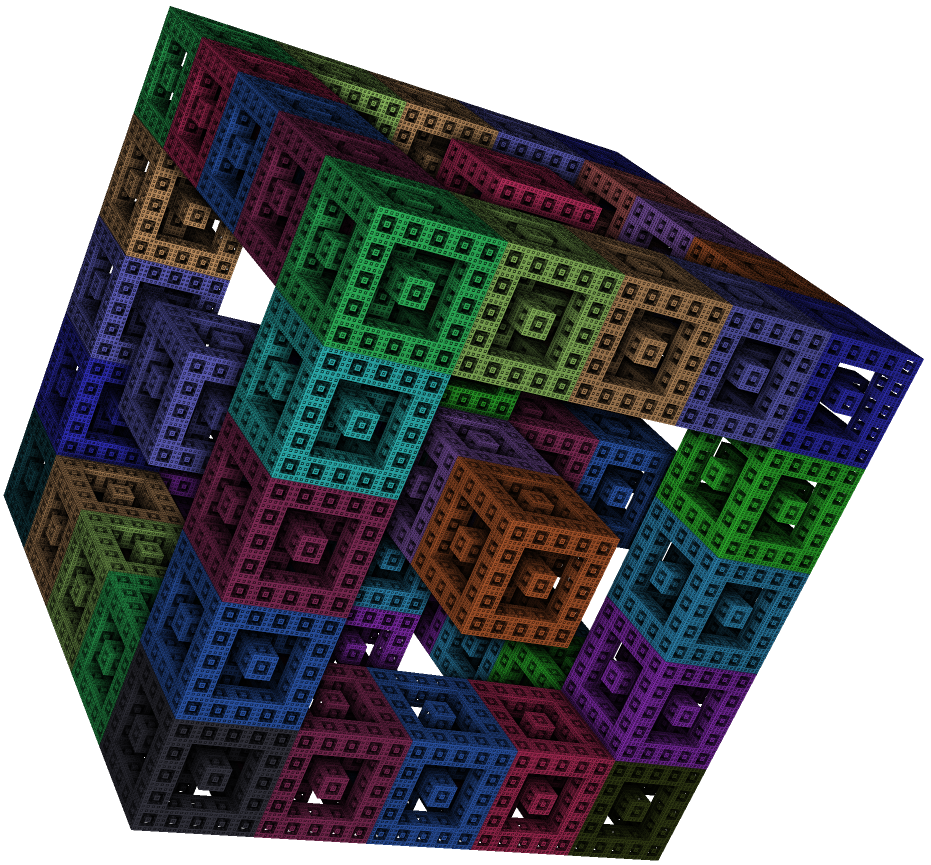}
\begin{center}
{\scriptsize The set $T(I)$ (left)  and the fractal cube $F$ (right).}    
\end{center}

We see from the construction that $\#\mD_0=13,\ \#\mD_1=44.$ Moreover, the Hausdorff distance between $\mD_0+I$ and $\mD_1+I$ is $2\sqrt{2},$ while the minimal distance $d(\mD_0+I,\mD_1+I)$ between the points of those sets is equal to 1.

\includegraphics[width=.45\textwidth]{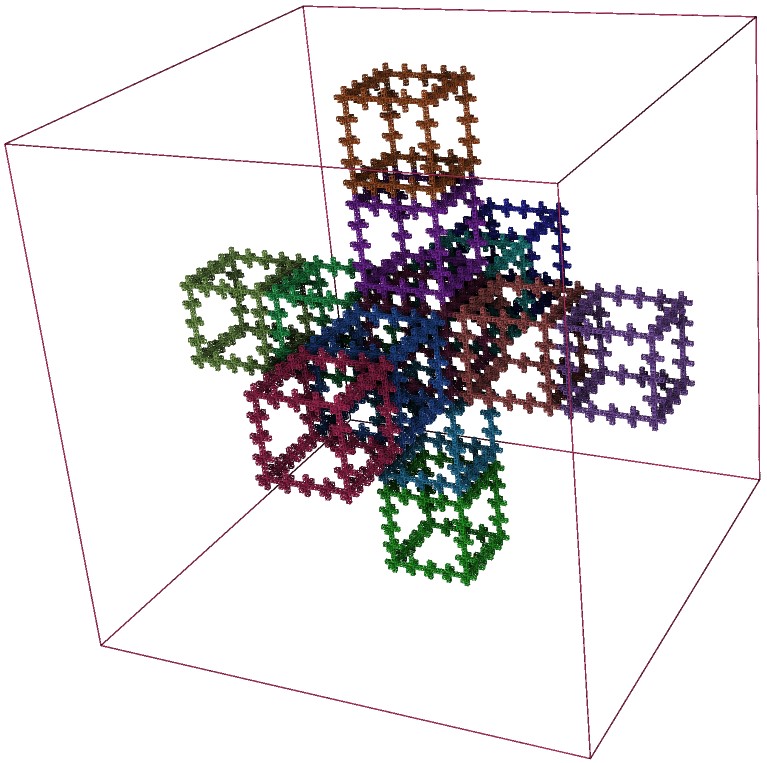}
\includegraphics[width=.45\textwidth]{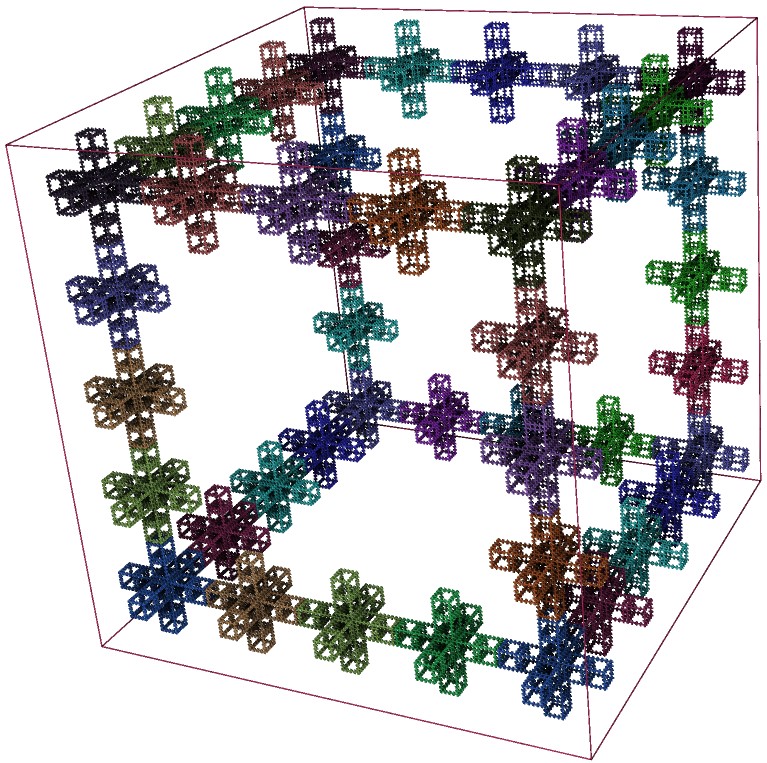}
\begin{center}
{\scriptsize The components $K_{01}$ and $K_{10}$}    
\end{center}
\begin{lem}
For any $\bar{\al},$ the sets  $F_{\bar{\al}}$ and $\eK_{\bar{\al}}$ are connected.
\end{lem}

\dok For $i=0,1,$ the intersections of $F_i$ with any two opposite faces of the cube $I$ are congruent with respect to the translation moving one face to the other.

For any two cubes $\dfrac{a+I}{5}, \dfrac{b+I}{5} \IN F_i$ with side $1/5$ in $F_i$ there is a path from one cube to the other which intersects transversely  the faces of pairs of adjacent cubes and does not intersect any of their edges.

It follows then, that for any $i,j$ the sets $T_i(F_j)$ are connected and possess the same opposite face property.

Proceeding by induction, we get  that the sets $T_{\al_1\ldots\al_k}(F_i)$ are connected and inherit the same property.

Since $\eK_{\bar{\al}}$ is the intersection of a nested family of compact connected sets $(T_{\bar{\al}})^m(I), \eK_{\bar{\al}}$ is connected too. \vse

\begin{lem}
If \quad $\bar{\al}, \bar{\be}\in\{0,1\}^k$ and $\bar{\al}\neq\bar{\be}$,  then\\     $F_{\bar{\al}}\bigcap F_{\bar{\be}}=\0,\ d_H(F_{\bar{\al}},F_{\bar{\be}})<\dfrac{3\sqrt5}{5^{|\bar{\al}\wedge\bar{\be}|}},$ and $d(F_{\bar{\al}},F_{\bar{\be}})\geq\dfrac{1}{5^{|\bar{\al}\wedge\bar{\be}|+1}}$
\end{lem}

\dok If $\al_1\neq\be_1$ then $F_{\bar{\al}}\IN F_{\al_1}, F_{\bar{\be}}\IN F_{\be_1},$ therefore 
\begin{equation}
    \begin{split}d_H(F_{\bar{\al}},F_{\bar{\be}})\leq d_H(F_{\bar{\al}},F_{\al_1})&+d_H(F_{\al_1},F_{\be_1})+d_H(F_{\bar{\be}},F_{\be_1})\\ & \leq \frac{2\sqrt2}{25}+\frac{2\sqrt2}{5}+\frac{2\sqrt2}{25} < \frac{3\sqrt2}{5}\end{split} \end{equation}
and $d(F_{\bar{\al}},F_{\bar{\be}})\geq d(F_{\al_1},F_{\be_1})>1/5.$ 
Notice, that $d_H(F_{\sa\al},F_{\sa\be})\leq\dfrac{d_H(F_{\al},F_{\be})}{5^{|\sa|}}$ to get the desired statement. \vse

Let $\bar{\al}=\al_1\al_2\al_3\ldots\in\{0,1\}^\8.$ We define $\eK_{\bar{\al}}=\bigcap\limits_{k=1}^\8 F_{\al_1\ldots\al_k}.$

\begin{cor}
For any $\bar{\al}, \bar{\be}\in\{0,1\}^\8$ 
$$\frac{1}{5^{|\bar{\al}\wedge\bar{\be}|+1}} \leq d_H(F_{\bar{\al}},F_{\bar{\be}}) < \frac{3\sqrt5}{5^{|\bar{\al}\wedge\bar{\be}|+1}},$$ and 
$$d(F_{\bar{\al}},F_{\bar{\be}})\geq \frac{1}{5^{|\bar{\al}\wedge\bar{\be}|+1}}$$
\end{cor}\vse


Consider $\bar{\al}\in\{0,1\}^k, k\in\nn$. The number $N_\da(F_{\bar{\al}})$ of $\da$-mesh cubes for $(F_{\bar{\al}}=(T_{\bar{\al}}(I)$ and $\da=5^{-k}$ is equal to $13^m\cdot 44^{k-m}$, where $m=\#\{i\leq k:\al_i=0\}$. Thus $-\dfrac{\log N_\da}{\log\da}=\dfrac{m}{k}\log_513+(1-\dfrac{m}{k})\log_544$.

For any $\al\in\{0,1\}^\8$, the minimal number of $\da$-mesh cube containing $F_\al$ is equal to $N_\da(F_{\bar{\al}})$ if $\da=5^{-k}$ and $\bar{\al}=\al_1\ldots\al_k$ is the $k$-th initial segment of $\al$. Taking upper and lower limits as $k\to\8$ we get desired estimates for $\dim_B(\eK_\al)$.

\end{document}